\documentclass[%
twocolumn
]{mpi2015-cscpreprint}

\usepackage{graphicx}      % include this line if your document contains figures

\usepackage[american]{babel}

\usepackage{amsmath, amssymb}
\usepackage{todonotes}

\usepackage{geometry}
\usepackage{verbatim}

\newtheorem{remark}{Remark}
\newcommand{\bA}{{\mathbf A}}
\newcommand{\bB}{{\mathbf B}}
\newcommand{\bC}{{\mathbf C}}

\newcommand{\bE}{{\mathbf E}}

\newcommand{\bG}{{\mathbf G}}
\newcommand{\bH}{{\mathbf H}}

\newcommand{\bQ}{{\mathbf Q}}

\newcommand{\bV}{{\mathbf V}}
\newcommand{\bW}{{\mathbf W}}
\newcommand{\bX}{{\mathbf X}}
\newcommand{\bY}{{\mathbf Y}}
\newcommand{\bJ}{{\mathbf J}}
\newcommand{\bK}{{\mathbf K}}

\newcommand{\bfe}{{\mathbf e}}

\newcommand{\bu}{{\mathbf u}}
\newcommand{\bv}{{\mathbf v}}
\newcommand{\bw}{{\mathbf w}}
\newcommand{\bx}{{\mathbf x}}
\newcommand{\by}{{\mathbf y}}

\newcommand{\cF}{{\mathcal F}}
\newcommand{\cT}{{\mathcal T}}
\newcommand{\cV}{{\mathcal V}}
\newcommand{\cZ}{{\mathcal Z}}
\newcommand{\bPhi}{ \boldsymbol{\Phi} }
\def\IR{{\mathbb R}}
\def\IC{{\mathbb C}}

\def\IL{{\mathbb L}}
\def\IV{{\mathbb V}}
\def\IW{{\mathbb W}}
\newcommand{\sIL}{{{{\mathbb L}_s}}}

\newtheorem{proposition}{Proposition}[section]

\begin{document}
	
	\title{Learning reduced-order models of quadratic control systems from input-output data} 
	
	\author[$\ast$]{Ion Victor Gosea}
	\affil[$\ast$]{Max Planck Institute for Dynamics of Complex Technical Systems, Magdeburg, Germany.\authorcr
		\email{gosea@mpi-magdeburg.mpg.de}, \orcid{0000-0003-3580-4116}}
	
	\author[$\dagger$]{Dimitrios S. Karachalios}
	\affil[$\dagger$]{Max Planck Institute for Dynamics of Complex Technical Systems, Magdeburg, Germany.\authorcr
		\email{karachalios@mpi-magdeburg.mpg.de}, \orcid{0000-0001-9566-0076}}
	
	\author[$\ddagger$]{Athanasios C. Antoulas}
	\affil[$\ddagger$]{Electrical and Computer Engineering (ECE) Department, Rice University, Houston, USA,
		
		 Max Planck Institute, Magdeburg, Germany, and Baylor College of Medicine, Houston, USA.\authorcr
		\email{aca@rice.edu}}
	
	\shorttitle{Learning quadratic models from data}
	\shortauthor{I.V. Gosea, D. S. Karachalios, and A. C. Antoulas}
	\shortdate{}
	

	\abstract{%
	In this paper, we address an extension of the Loewner framework for learning quadratic control systems from input-output data. The proposed method first constructs a reduced-order linear model from measurements of the classical transfer function. Then, this surrogate model is enhanced by incorporating a term that depends quadratically on the state. More precisely, we employ an iterative procedure based on least squares fitting that takes into account measured or computed data. Here, data represent transfer function values inferred from  higher harmonics of the observed output, when the control input is purely oscillatory. }

		\maketitle

%\begin{keyword}
%Data-driven modeling, nonlinear systems,  input-output data, reduced-order modeling, Loewner matrix, transfer functions, model reduction, frequency-domain methods.
%\end{keyword}

%\end{frontmatter}
%===============================================================================

\section{Introduction}

With an ever-increasing availability of measured data in many engineering fields, the need for incorporating measurements in the modeling process has steadily grown over the last decades. The main challenge lies within effectively using the available data in order to construct models that can accurately represent the dynamics of the underlying dynamical process. Sometimes, in order to satisfy accuracy requirements, the fitted models might have large dimension and hence are not suitable for fast numerical simulation. Hence, it is of interest to compute reliable reduced-order surrogate models instead. 

Model reduction is commonly viewed as a
methodology used for reducing the computational complexity of large scale complex models in
numerical simulations. The goal is to construct a smaller system with the same structure
and similar response characteristics as the original. For an overview of conventional model reduction methods, we refer the reader to the books \cite{ACA05,BOCW17,ABG20}. In this work, we assume that the nonlinear systems to be modeled contains quadratic nonlinearities. This class of systems is of interest since most smooth nonlinear systems can be exactly reformulated as quadratic or quadratic-bilinear (QB) systems (provided that the nonlinearities are analytical). MOR methods specifically tailored for reducing QB systems have been proposed in \cite{BB15,BGG18}. For a general overview on system-theoretical nonlinear MOR approaches, see \cite{BBF14}.

This work concentrates on non-intrusive model reduction methods. These represent a class of methodologies that do not necessarily require access to the original model (described, e.g., by matrices or other operators), but only to measured data (in the frequency or in time domain). The method that is in the center of the current study is the Loewner Framework (LF). It is a data-driven model identification and reduction technique that was originally introduced in \cite{MA07}. Using only measured input-output data in the frequency domain, it directly constructs reliable surrogate models. Other frequency-domain methods include vector fitting in \cite{GS99} or the moment matching method in \cite{astolfi2010}.

Non-intrusive methods that use time-domain data include, e.g., subspace identification in \cite{OM94,Qin06}.
Recently, other methods emerged, such as dynamic mode decomposition (DMD) in \cite{morSch10,PBK16}, and operator inference (OpInf) in \cite{PEHERSTORFER2016196,morBenGKetal20}. The latter is used to learn models by means of fitting reduced-order quadratic operators that minimize a time-domain deviation in a least squares (LS) sense. It requires access only to snapshots of the states and its derivatives.

The method proposed in this work comes as a continuation of \cite{morGosKA20}, and represents an extension of the case with bilinear nonlinearities treated in \cite{morKarGA20}. It is a non-intrusive method in the sense that we do not require access to the system's matrices. Only measurements of higher-order transfer functions are needed (harmonic content of the observed output). To fit the linear part, we use first transfer function values only. Information from the second and the third harmonics is used to fit the quadratic part that best complements the already-fitted linear model. At first, this is done by means of a one-step approach that fits the harmonics separately. Afterwards, we propose an iterative algorithm that solves a linearized LS coupled problem (by taking into account all data). Different than earlier methods such as \cite{morGosKA20,morKarGA20}, in this paper we propose an adaptive iterative scheme for learning reduced-order models based on fitting input-output data corresponding to higher-order harmonics. We also take into account noisy data, in the sense that the measurements can be assumed to be corrupted by noise. 

The paper is organized as follows; after the introduction, a section on quadratic control systems is provided in Section \ref{sec:2}. It presents the state-space description of these systems together with the definition of its generalized transfer functions. Then, the classical Loewner framework for modeling linear systems is introduced in Section \ref{sec:loew}. The proposed method is detailed in Section \ref{sec:recQ} and summarized in Algorithm 1. In Section \ref{sec:5} we present numerical results for applying the new method on two test cases (low and high-dimensional). Section \ref{sec:6} states the conclusions and potential future developments, while Section \ref{sec:7} includes an appendix.

\section{Quadratic control systems}
\label{sec:2}

\noindent
Consider quadratic control systems that are characterized by the following equations
\begin{align}\label{quad_sys}
\begin{split}
\bE\dot\bx(t)&=\bA\bx(t)+\bQ\,(\bx(t)\otimes\bx(t))+\bB\bu(t),\\
\by(t)&=\bC\bx(t).
\end{split}
\end{align}
with $\bE, \bA \in \IR^{n \times n}, \bB, \bC^T \in \IR^{n}, \bQ \in \IR^{n \times n^2}$, $u(t)$ is the control input and $y(t)$ is the observed output. The Kronecker product of the internal variable $\bx = [x_1 \  \cdots \  x_n]^T$ with itself is used in (\ref{quad_sys}), i.e.  
$$
\bx \otimes \bx = [ x_1^2 \ \ x_1 x_2  \ \ \ldots \ \ x_1 x_n 
\ \ \ldots \ \  x_n^2]^T \in \IR^{n^2}.
$$ 
Note also that the matrix $\bE$ is considered to be invertible. For such class of systems, one can explicitly compute generalized transfer functions in the frequency domain. This is done based on Volterra series theory (we refer the reader to \cite{WJR05} for details). In general, the Volterra series describes the relationship between the control input and the observed output of a dynamical system with nonlinear dynamics. An explicit formulation of input-output mappings in the time domain is provided by means of the Volterra kernels. The frequency domain equivalent of these kernels is represented by the generalized transfer functions (which are multi-variate rational functions). Such functions can be derived from time-domain data by applying generalized Fourier or Laplace transformations.

Deriving analytical expressions of generalized transfer functions can be made using the harmonic input probing method in \cite{BR71}. It is based on the fact that a harmonic input must result in a harmonic output. More details on approximating nonlinear operators by Volterra series are given in \cite{BoydChua85}.\\
For simplicity, in what follows, we use a single-tone harmonic function  to identify explicit formulas for the generalized transfer functions in the case of quadratic systems.

%Consequently, the transfer function be evaluated at multiple instances of the driving frequency $j \omega$.

The input signal is purely oscillatory, i.e.
$\bu(t) = \alpha e^{j \omega t}$, where $\omega, \alpha >0$ and $j = \sqrt{-1}$.
%We begin by
%making the following classical assumption for the solution of (\ref{quad_sys})
We start by assuming that the solution of (\ref{quad_sys}), i.e. $\bx(t)$, can be expanded in the following power series (as done, e.g., in \cite{WJR05})
\begin{equation} \label{sol_quad}
\bx(t) = \sum_{m=1}^\infty \bG_m(j\omega) \alpha^
m e^{m j \omega t},
\end{equation}
where $\bG_m: \mathbb{C} \rightarrow \mathbb{C}^n$. For an analysis on convergence for series as in (\ref{sol_quad}), we refer the reader to \cite{WJR05}.

Substitute the relation (\ref{sol_quad}) into the original differential equation (\ref{quad_sys}), and hence write that:
\begin{align}\label{quad_diff_eq}
&\bE \sum_{m=1}^\infty m j \omega \bG_m(j\omega) \alpha^m e^{m j \omega t} = \bA \sum_{m=1}^\infty \bG_m(j\omega) \alpha^m e^{m j \omega t} \nonumber \\
&+\bQ \Big{(} \sum_{m=1}^\infty \bG_m(j\omega) \alpha^m e^{m j \omega t} \Big{)}  \otimes \Big{(} \sum_{m=1}^\infty \bG_m(j\omega) \alpha^m e^{m j \omega t} \Big{)} \nonumber \\ &+ \bB \alpha e^{j \omega t}. 
\end{align}
Next, in order to explicitly find a closed-form expression for $\bG_m(j\omega)$, one needs to equate the coefficient of the term $e^{m j \omega t}, \ \forall m \geqslant 1$ from both left and right sides in (\ref{quad_diff_eq}).\\
For $m=1$, it follows that $\bG_1(j \omega ) = \bPhi(j\omega) \bB$ is identified, where $\bPhi(j\omega) = (j \omega \bE - \bA)^{-1}$. The first transfer function is hence given by
\begin{align}\label{first_trf}
\bH_1(j \omega ) = 
\bC \bG_1(j \omega ) = \bC \bPhi(j\omega) \bB.
\end{align}
Next, for $m=2$, one can write that
\begin{align*}
\begin{split}
& 2j \omega \bE \bG_2(j \omega) \alpha^2  = \bA \bG_2(j \omega) \alpha^2 + \bQ [ \bG_1(j \omega) \otimes \bG_1(j \omega)]  \alpha^2\\
& \Rightarrow  \bG_2(j \omega)  = \bPhi(2j\omega) \bQ [\bG_1(j \omega) \otimes \bG_1(j \omega)].
\end{split}
\end{align*}
Hence, write the second transfer function as
\begin{equation} \label{second_trf_quad}
\bH_2(j \omega) = \bC \bPhi(2j\omega) \bQ \left[\bPhi(j\omega) \bB \otimes  \bPhi(j\omega) \bB \right].
\end{equation}
By repeating this procedure, for $m=3$, one can write that
%\small
%\begin{align*}
%\begin{split}
%  & \bG_3(j \omega) = \bPhi(3j\omega) \bQ [\bG_2(j\omega) \bB \otimes \bG_1 (j \omega) + \bG_1 (j \omega) \otimes \bG_2(j\omega)].
%\end{split}
%\end{align*}
%\normalsize
%Hence, the third transfer function is written as follows
\begin{align} \label{third_trf_quad}
\begin{split}
\bH_3(j \omega) &= \bC \bPhi(3j\omega) \bQ \left[\bG_2(j\omega) \otimes \bG_1 (j \omega)\right] \\ &+ \bC \bPhi(3j\omega) \bQ \left[ \bG_1 (j \omega) \otimes \bG_2(j\omega)\right].
\end{split}
\end{align}
By imposing the symmetry rule $\bQ(\bv \otimes \bw) = \bQ(\bw \otimes \bv)$, one can simplify the formula above as follows, i.e.,
\begin{align} \label{third_trf_quad2}
\begin{split}
& \bH_3(j \omega) = 2\bC \bPhi(3j\omega) \bQ \left[\bG_2(j\omega) \otimes \bG_1 (j \omega)\right]
%&= \bC \bPhi(3j\omega) \bQ \left[\bPhi(j\omega) \bB \otimes  \bPhi(2j\omega) \right] \bQ \left[\bPhi(j\omega) \bB \otimes  \bPhi(j\omega) \bB \right]
\end{split}
\end{align}
In general, the $n$th transfer function is written as $\bH_m(j \omega) = \bC \bG_m(j\omega)$ and can be expressed in closed-form for any $n \geq 1$,  following a similar procedure as for that used to derive (\ref{second_trf_quad}). Based on (\ref{sol_quad}), one can write the time-domain observed output in (\ref{quad_sys}) as follows
\begin{equation} \label{output_quad_time}
y(t) = \bC \bx(t) =  \sum_{m=1}^\infty \bH_m(j\omega) \alpha^
m e^{m j \omega t},
\end{equation}
while the frequency-domain representation is obtained by applying the Fourier transform ($\cF(\cdot)$) to y(t), as
\begin{equation} \label{output_quad_freq}
Y(j \Omega) = \cF(y(t)) =  \sum_{m=1}^\infty \bH_m(j\omega) \alpha^
m \delta(j(\Omega-\omega)),
\end{equation}
where $\delta(\cdot)$ is the Dirac delta function.
Hence, by analyzing the frequency content of the observed output explicitly given in (\ref{output_quad_freq}), one can estimate transfer function values $\bH_m(j\omega)$ by means of the harmonic content. More precisely, we use that $\bH_m(j \omega) = \cZ_m / \alpha^m$, where $\cZ_m$ is the mth harmonic.

\section{The Loewner framework}\label{sec:loew}

In this section we present a short summary of the Loenwer framework (LF). For a tutorial paper on LF for linear systems, we refer the reader to \cite{ALI17}. For an extension that uses time-domain data, see \cite{PGW17}. The Loewner framework has been recently extended to certain classes of nonlinear systems, such as bilinear systems in \cite{AGI16}, and quadratic-bilinear (QB) systems in \cite{GA18}.

The starting point for the classical LF is to collect measurements corresponding to the (first) transfer function, which can be inferred in practice from the first harmonic. The problem is formulated as follows.
Given the following scalar data values partitioned into two disjoint subsets
\begin{align}\label{data_Loew}
\begin{split}
{\text{ right \ data}}: &(\lambda_i;\bw_i), ~i=1,\ldots,k,~~{\text and}, \\
{\text{ left \ data}}: &(\mu_j;\bv_j), ~j=1,\ldots,k,
\end{split}
\end{align} 
%(for simplicity all points are assumed distinct).
find the function
$\bH(s)$, such that the following interpolation conditions are (approximately) fulfilled:
\begin{equation} \label{interp_cond}
\bH(\lambda_i)=\bw_i,~~~\bH(\mu_j)=\bv_j.
\end{equation}
The Loewner matrix $\IL \in\IC^{k\times k}$ and the shifted Loewner matrix $\sIL \in\IC^{k\times k}$ are defined as follows
\begin{equation} \label{Loew_mat}
\IL_{(i,j)}=\frac{\bv_i-\bw_j}{\mu_i-\lambda_j}, \ \sIL_{(i,j)}=
\frac{\mu_i\bv_i-\lambda_j \bw_j}{\mu_i-\lambda_j},
\end{equation}
while the data vectors $\bV, \bW^T \in \IR^k$ are introduced as
\begin{equation} \label{VW_vec}
\IV_{(i)}= \bv_i, \ \  \IW_{(j)} = \bw_j.
\end{equation}

The Loewner model is composed of
\begin{align*}
\bE=-\IL,~~ \bA=-\sIL,~~ \bB=\IV,~~ \bC=\IW.
\end{align*}

In practical applications, the pencil $(\sIL,\,\IL)$ is often singular. In these cases, perform a rank revealing singular value decomposition (SVD) of the Loewner matrices. Then, compute projection matrices $\bX_r, \bY_r \in \IC^{k \times r}$ as the left, and respectively, the right truncated singular vector matrices. More details can be found in \cite{ALI17}. Here, $r<n$ represents the truncation index.
%\begin{equation*}
%~\IL =\bX \bS \bY^* \approx \bX_r\bS_r \bY_r^*, \ \ \ %\text{with} \ \ \bX_r, \bY_r \in \IC^{k \times r}, \  %\bS_r \in \IC^{r \times r}.
%\end{equation*}	
Then, the system matrices corresponding to a projected Loewner model of dimension $r$ can be computed
\begin{equation}\label{Loew_red_lin}
\hat{\bE} = -\bX_r^*\IL \bY_r, \ \  \hat{\bA} = -\bX_r^*\sIL \bY_r, \ \
\hat{\bB} = \bX_r^*\IV, \ \  \hat{\bC} = \IW \bY_r.
\end{equation}

\section{The proposed method for learning a quadratic model from data}\label{sec:recQ}

% In this contribution, the proposed procedure represents a natural extension of the method introduced in \cite{KGA19} from the case of bilinear systems to the case of quadratic systems.

In this section we present the main method. As mentioned before, the contribution is based on multiple-steps learning. After going through several separated stages of approximation, we propose an algorithm for solving a coupled minimization problem (that takes into account both second and third harmonic content). The main difference between the newly-proposed method and that in \cite{GA18}, is that the former takes into account measurements that can be inferred from experiments, i.e., from the frequency content of the observed output.

\subsection{First step: learning a linear surrogate model}\label{sec:4.1}

Using the classical Loewner framework introduced in Section \ref{sec:loew}, we first fit a linear model that matches samples of the first (linear) transfer function. Then, from samples of the second transfer function (that includes the nonlinear behavior), we are able to fit appropriate quadratic terms.
%The data required for this procedure can be estimated from direct numerical simulations in the time domain.

One can use the classical Loewner framework approach to directly construct a reliable reduced-order linear model $(\tilde{\bA},\tilde{\bB},\tilde{\bC})$ of order $r$ from samples of the first transfer function $\bH_1(j \omega)$ in (\ref{first_trf}), evaluated at the $N$ sampling points  $\{j\omega_1,\ldots,j \omega_{N}\}$.  The way to compute these new matrices is by simplifying the description in (\ref{Loew_red_lin}), as follows
\begin{equation}\label{Loew_red_lin2}
\tilde{\bA} = \hat{\bE}^{-1} \hat{\bA}\ \
\tilde{\bB} = \hat{\bE}^{-1} \hat{\bB}, \ \  \tilde{\bC} = \hat{\bC}.
\end{equation}

%These sample values are stored in the vector $\cV^{(1)} \in \IC^{N}$.

For ease of notation, in the next sections we refer to the projected Loewner model as to the (linear) reference model and use the notation $(\bA,\bB,\bC)$ instead of the one in (\ref{Loew_red_lin2}).

\subsection{Learning a reduced-order quadratic operator from samples of the second transfer function}
\label{sec:4.2}

The next step is to fit an appropriate matrix $\bQ \in \IC^{r \times r^2}$ that supplements the linear model into a quadratic model. In this direction, it is assumed that information about the second transfer function in (\ref{second_trf_quad}) is known at $N$ points $\{j\omega_1,\ldots,j \omega_{N}\}$. These measured sample values are stored in the vector $\cV^{(2)} \in \IC^{N}$.

\begin{remark}
	Higher-order transfer function values can be estimated from transformed output trajectories, when the control input is a purely oscillating signal. In this case, we need to choose N different control inputs $u_m(t) = e^{j \omega_m t}$ for $1 \leq m \leq N$, and then simulate the black-box model.
\end{remark}

The problem is to find a reduced-order matrix $\bQ  \in \IC^{r \times r^2}$ so that to minimize the distance between the data corresponding to the original model and the data corresponding to the learned low-dimensional model.
% \begin{equation}
% \sum_{\ell=1}^{N} \left( \bH_2(j\omega_\ell) - \cV^{(2)}_\ell \right)^2.
% \end{equation}
We hence set up the minimization problem in $r^3$ unknowns
\begin{equation}\label{min_prob_2TF}
\min_{\bQ \in \IC^{r \times r^2}}  \sum_{\ell=1}^{N}  \left( \bH_2(j\omega_\ell) - \cV^{(2)}_\ell \right)^2.
\end{equation}
Let $\bJ_1: \mathbb{C} \rightarrow \mathbb{C}^{1 \times r}$, such that $\bJ_1(j \omega) = \bC \bPhi(j\omega)$. Denote with $\bv_\bQ \in \IC^{r^3}$ the vectorization of the matrix $\bQ  \in \IC^{r \times r^2}$:
\begin{equation} \label{vecQ}
\bv_\bQ = \left[(\bQ \bfe_1)^T \  (\bQ \bfe_2)^T \ \cdots \  \ (\bQ \bfe_{r^2})^T\right]^T.
\end{equation}
Here, $\bfe_h \in \mathbb{C}^{r^2}$ is the $h$th unit vector. 
% Then, the following holds true
\begin{proposition}\label{prop:1}
	The minimization problem in (\ref{min_prob_2TF}) can be equivalently written as
	\begin{equation}\label{min_prob_2TF_2}
	\min_{\bv_\bQ \in \IC^{r^3}} \Vert \mathcal{T}^{(2)} \bv_\bQ - \mathcal{V}^{(2)} \Vert^2,
	\end{equation}
	%\begin{equation} \label{eq_vQ}
	%\mathcal{T}^{(2)} \bv_\bQ = \mathcal{V}_2,  
	%\end{equation}
\end{proposition}
where $\cT^{(2)} \in \IC^{N \times r^3}$ with its $\ell$th row given by ($1 \leq \ell \leq N$)
\begin{equation}\label{def_T2}
\bfe_\ell^T \cT^{(2)} =  \bG_1^T(j \omega_\ell) \otimes \bG_1^T(j \omega_\ell) \otimes \bJ_1(2j \omega_\ell).
\end{equation}
Note that the matrix $\cT^{(2)}$ introduced in (\ref{def_T2}) depends entirely on the linear (Loewner) model constructed in Section \ref{sec:4.1}.

Hence, the problem in (\ref{min_prob_2TF_2}) is linear in vector $\bv_\bQ \in \IC^{r^3}$, and can be directly solved by means of, for example, the Moore-Penrose pseudo-inverse matrix ${\cT^{(2)}}^{\dagger}$. More precisely, we can write the solution vector as
\begin{equation}\label{estimateQ2}
\bv_\bQ^{(2)} = {\cT^{(2)}}^{\dagger} \cV^{(2)}.
\end{equation}
Afterwards, one can automatically reconstruct the recovered matrix $\bQ \in \IC^{r \times r^2}$ based on the formula in (\ref{vecQ}).
\begin{remark}
	We are using additional harmonic content instead of relying on the model learned at this stage, given by (\ref{estimateQ2}). The motivation is that this might not provide enough information to accurately represent the underlying dynamics. 
	%More precisely, matching the first two TFs could be achieved by a broader family of surrogate models.
\end{remark}

\subsection{A minimization problem formulated using samples of the third transfer function} \label{sec:4.3}

In this cases, it is assumed that values of the third transfer function in (\ref{third_trf_quad}) are estimated (at the same points as before, i.e., $\{j\omega_1,\ldots,j \omega_{N}\}$). These measured sample values are stored in the vector $\cV^{(3)} \in \IC^{N}$. 
%The problem here consists in finding a reduced-order quadratic operator $\bQ  \in \IC^{r \times r^2}$ that minimizes the distance between the data corresponding to the original model and the data corresponding to the learned low-dimensional model:   
The problem here is formulated similarly as to that in Section \ref{sec:4.2}, but instead using measurements of the third transfer function, i.e.,
\begin{equation}\label{min_prob_3TF}
\sum_{\ell=1}^{N} \left( \bH_3(j\omega_\ell) - \cV^{(3)}_\ell \right)^2,
\end{equation}
Let $\bK: \mathbb{C} \rightarrow \mathbb{C}^{r^3 \times r^3}$, explicitly given by the formula:
\begin{align}\label{def_K}
\begin{split}
\bK(j\omega) &= 
2\left[ \left( \bG_1(j \omega) \otimes \bG_1(j \omega) \right) \otimes \bPhi^T(2j \omega) \right] \\
& \otimes \bG_1^T(j \omega) \otimes \bJ_1(3j \omega). 
\end{split}
\end{align}
\begin{proposition}\label{prop:2}
	Minimizing the quantity in (\ref{min_prob_3TF}) is equivalently formulated as follows:
	\begin{equation}\label{min_prob_3TF_2}
	\min_{\bv_\bQ \in \IC^{r^3}}  \sum_{\ell=1}^{N} \left(  \bv_{\bQ}^T  \bK(j\omega_\ell)  \bv_\bQ - \cV^{(3)}_\ell \right)^2.
	\end{equation}
\end{proposition}
\vspace{3mm}
Note that the mapping $\bK$ introduced in (\ref{def_K}) depends entirely on the linear (Loewner) model constructed in Section \ref{sec:4.1}. The problem in (\ref{min_prob_3TF_2}) is hence quadratic in the variable $\bv_\bQ \in \IC^{r^3}$. To linearize it, one can use the estimate of $\bv_\bQ$ in (\ref{estimateQ2}) and replace it in (\ref{min_prob_3TF_2}), i.e.,
\begin{align}\label{min_prob_3TF_lin}
\begin{split}
\min_{\bv_\bQ \in \IC^{r^3}}  \sum_{\ell=1}^{N} \left(  \left(\bv_{\bQ}^{(2)}\right)^T  \bK(j\omega_\ell)  \bv_\bQ - \cV^{(3)}_\ell \right)^2 \\ \Leftrightarrow   \min_{\bv_\bQ \in \IC^{r^3}} \Vert \mathcal{T}^{(3)} \bv_\bQ - \mathcal{V}^{(3)} \Vert^2.
\end{split}
\end{align}
In (\ref{min_prob_3TF_lin}), the matrix $\cT^{(3)} \in \IC^{N \times r^3}$ is introduced, with its $\ell$th row given by $\left(\bv_{\bQ}^{(2)}\right)^T  \bK(j\omega_\ell)$.
As before, the problem in (\ref{min_prob_3TF_lin}) is linear in vector $\bv_\bQ \in \IC^{r^3}$, and can be directly solved by using the Moore-Penrose inverse denoted with ${\cT^{(3)}}^{\dagger}$, as
\begin{equation}\label{estimateQ3}
\bv_\bQ^{(3)} = {\cT^{(3)}}^{\dagger} \cV_3.
\end{equation}

\subsection{Solving the coupled problem by means of an iteration} \label{sec:4.4}

In this section we aim at learning a reduced-order quadratic operator $\bQ  \in \IC^{r \times r^2}$ that minimizes the coupled distance between data corresponding to second and third  transfer functions, as described below   
\vspace{-1mm}
\begin{equation}
\sum_{\ell=1}^{N} \left( \bH_2(j\omega_\ell) - \cV^{(2)}_\ell \right)^2 +  \sum_{\ell=1}^{N} \left( \bH_3(j\omega_\ell) - \cV^{(3)}_\ell \right)^2. 
\vspace{-1mm}
\end{equation}
Based on the results presented in the previous two sections, the coupled minimization problem is hence written as
\vspace{-2mm}
\begin{align}\label{min_prob_coupl}
\begin{split}
&\min_{\bv_\bQ \in \IC^{r^3}}  \Vert \mathcal{T}^{(2)} \bv_\bQ - \mathcal{V}^{(2)} \Vert^2 + \sum_{\ell=1}^{N} \left(  \bv_{\bQ}^T  \bK(j\omega_\ell)  \bv_\bQ - \cV^{(3)}_\ell \right)^2.
\end{split}
\end{align}
Next, we present an algorithm meant to solve the nonlinear problem stated in (\ref{min_prob_coupl}), by means of a fixed-point iteration scheme. Let $\tau >0$ be a tolerance value.

%\newpage

\vspace{2mm}

\textbf{Algorithm 1}
\vspace{2mm}
\begin{enumerate}
	\item Initialization step ($q=0$): first solve the linear part in (\ref{min_prob_coupl}), i.e. find the first estimate for $\bv_\bQ \in \IC^{r^3}$, as
	$\overline{\bv}_\bQ = {\cT^{(2)}}^{\dagger} \cV_2$.
	Also, set $\widetilde{\bv}_{\bQ} = 0$.\\
	
	\vspace{0mm}
	\noindent
	\textbf{While} $\vert  \overline{\bv}_\bQ -  \widetilde{\bv}_\bQ \vert > \tau$ \textbf{do} \\
	\item If $q \geqslant 1$, let $   \overline{\bv}_\bQ = \widetilde{\bv}_\bQ$ (update vector $ \overline{\bv}_\bQ$ with $ \widetilde{\bv}_\bQ$).\\
	\item Construct the matrix $\overline{\cT}^{(3)} \in \IC^{N \times r^3}$  such that its $\ell$th row is given by the following formula (for $1 \leq \ell \leq N$)
	\begin{equation}\label{def_T3_over}
	\bfe_\ell^T \overline{\cT}^{(3)} =  \left(  \overline{\bv}_\bQ\right)^T  \bK(j\omega_\ell),
	\end{equation}
	where $\bK: \mathbb{C} \rightarrow \mathbb{C}^{r^3 \times r^3}$ is as in (\ref{def_K}).
	\item The problem in (\ref{min_prob_coupl}) becomes linear; rewrite it as
	\begin{equation}\label{min_prob_coupl_lin}
	\min_{\bv_\bQ \in \IC^{r^3}} \Vert \mathcal{T}^{(2)} \bv_\bQ - \mathcal{V}^{(2)} \Vert^2+ \Vert \overline{\cT}^{(3)} \bv_\bQ - \mathcal{V}^{(3)} \Vert^2.
	\end{equation}
	\item Explicitly write the solution to (\ref{min_prob_coupl_lin}) as follows:
	\begin{equation}\label{tilde_v_Q}
	\widetilde{\bv}_\bQ = \left( \left[ \begin{matrix}
	\mathcal{T}^{(2)} \\ \overline{\cT}^{(3)}
	\end{matrix} \right] \right)^{\dagger} \left[ \begin{matrix}
	\mathcal{V}^{(2)} \\ \mathcal{V}^{(3)} 
	\end{matrix} \right].
	\end{equation}
	%    \item Repeat until the deviation between the $\widetilde{\bv}_\bQ$ and $\overline{\bv}_\bQ$ values drops below the tolerance value $\tau$.
	\item $q = q+1$.\\
	\textbf{end}
\end{enumerate}

%\vspace{2mm}

The output of this algorithm, upon convergence, is a vector ${\bv}_\bQ \in \IC^{r^3}$. Reshape it as a matrix $\bQ \in \IC^{r \times r^2}$ and then augment the Loewner linear model in Section \ref{sec:4.1} by a learned quadratic model $(\bA,\bB,\bC,\bQ)$.

\begin{remark}
	We compute pseudo-inverses throughout Section  \ref{sec:recQ} and also in Algorithm 1 be means of a SVD approach. More specifically, a threshold value $\epsilon>0$ is used to determine the singular values/vectors that actually enter the computations (similar to the OpInf procedure).
\end{remark}

\begin{remark}
	We assume that the transfer function measurements are known at evaluation frequencies $j \omega_\ell$ as well as $-j \omega_\ell$. By involving data closed under complex conjugation, we are able to compute real-valued surrogate models (as described for the linear case in \cite{ALI17}).
\end{remark}

\section{Test cases}\label{sec:5}

\vspace{-1mm}

\subsection{A low-dimensional example}

Consider a simple quadratic model, as in (\ref{quad_sys}), described by matrices
\begin{align}
\bA_o = \left[ \begin{matrix}
-0.03 & -2 \\ 2 & -0.05
\end{matrix} \right], \ \  \bQ_o = \left[ \begin{matrix}
1 & 0 & 0 & 0 \\ 0 & 0.5 & 0.5 & 0
\end{matrix} \right],
\end{align}
and by $\bB_o = \left[ \begin{matrix}
1 & 1
\end{matrix} \right]^T,  \bC_o = \left[ \begin{matrix}
1 & 0
\end{matrix} \right]$. The data consist of 40 logarithmically-spaced points (and complex conjugates) inside the interval $[10^{-0.5},5]$ together with evaluations of the first three transfer functions at these points. In this example, we treat the noise-free case. We will be using quadratic terms computed based on two approaches
\vspace{-1mm}
\begin{enumerate}
	\item the one-step solution provided in Section \ref{sec:4.2}.
	\item the solution of the iterative algorithm in Section \ref{sec:4.4}.
\end{enumerate}

\begin{figure}[h]
	\hspace{-6mm}
	\includegraphics[scale=0.245]{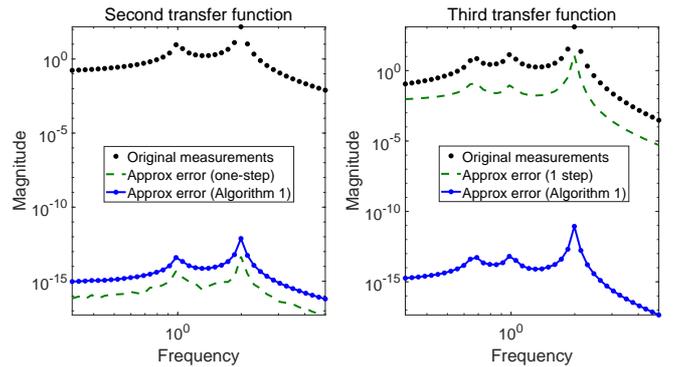}
	\vspace{-4mm}
	\caption{Original measurements and approximation errors.}
	\label{fig:1}
	\vspace{-2mm}
\end{figure}

Note that the linear part is perfectly recovered in the Loewner modeling stage (for both methods). The  one-step estimate computed as in (\ref{estimateQ2}) is given by:
\begin{equation}
\bQ^{(2)} =  
\left[ \begin{matrix}
0.9704  &  0.0854  &  0.0854  & -0.1412 \\
0.0594  &  0.3877  &  0.3877 &  -0.2789
\end{matrix} \right].
\end{equation}
The tolerance value is $\tau = 10^{-14}$ and after performing 45 steps of Algorithm 1, the deviation $\vert  \overline{\bv}_\bQ -  \widetilde{\bv}_\bQ \vert$ drops below the value $\tau$. We recover the matrix $\bQ$, with an error given by $\Vert \bQ - \bQ_o \Vert_2 = 2.0024 \cdot 10^{-14}$. Next, in Fig.\;\ref{fig:1} we depict the original measurements and the approximation errors of the second and of the third transfer functions corresponding to the two methods.

Note that the method in Algorithm1 recovers the original system, while the other one fails to approximate the third transfer function (as it can be seen in Fig.\;\ref{fig:1}, right pane).  

%\begin{figure}[h]		
%	\hspace{-6mm}		\includegraphics[scale=0.24]{plot_toy2.eps}
%	\vspace{-3mm}
%	\caption{....}
%	\label{fig:2}
%	\vspace{-2mm}
%\end{figure}

\vspace{-1mm}

\subsection{A large-scale MOR example}  
Consider a semi-discretized model of the Burgers' equation of dimension $n=100$ as described in \cite{BB15} (the bilinear part is neglected). Consider $100$ log-spaced points in the interval $2\pi[10^{-2},10^2]i$. Data consist in samples of the first three transfer functions. Note that the measurements are corrupted with Gaussian noise with signal to noise ratio (snr) value of 60dB.
In what follows, apply the method proposed in Algorithm 1.

\begin{figure}[h]
	\hspace{-6mm}
	\includegraphics[scale=0.245]{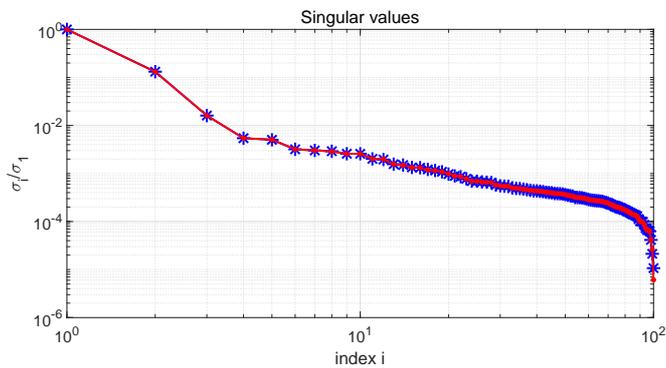}	\vspace{-4mm}
	\caption{The singular values of the Loewner matrix.}
	\label{fig:2}
	\vspace{-3mm}
\end{figure}

The reduction order is chosen to be $r=5$.
This choice is made with respect to the first neglected normalized singular value, i.e.,  $\sigma_{6}/\sigma_{1}\approx 10^{-3}$, as it can be seen in the decay of dingular values depicted in Fig.\;\ref{fig:2}.

\begin{figure}[h]
	\hspace{-6mm}		\includegraphics[scale=0.245]{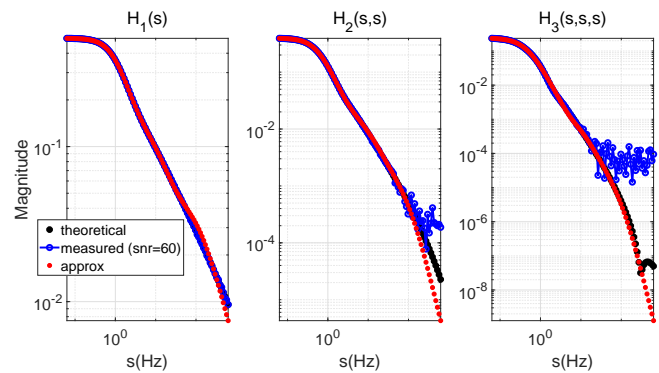}	\vspace{-4mm}
	\caption{The theoretical, measured and approximated values corresponding to the first  three transfer functions.}
	\label{fig:3}
	\vspace{-3mm}
\end{figure}

Since the LS problems could be under-determined, we are using a thresholding singular value decomposition scheme by choosing a threshold value of $\epsilon =  10^{-3}$. Note that this is chosen in accordance to the first neglected normalized singular of the Loewner matrix. In Fig.\;\ref{fig:3}, the three theoretical responses along with the noisy estimations, and also the transfer functions of the fitted reduced-model evaluated over the frequency grid.

\begin{figure}[h]
	\hspace{-6mm}		\includegraphics[scale=0.245]{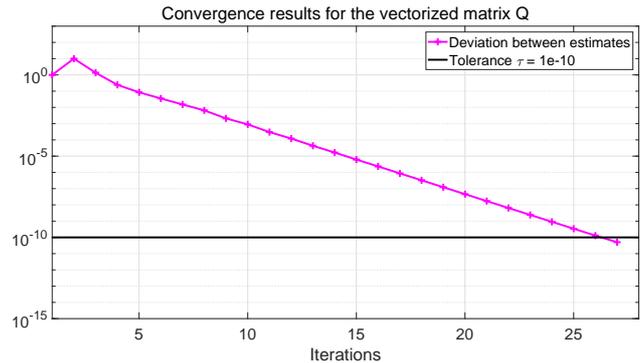}
	\vspace{-4mm}
	\caption{Convergence rate of Algorithm 1.}
	\label{fig:4}
	\vspace{-3mm}
\end{figure}

In Fig.\;\ref{fig:4} the exponential convergence rate is illustrated. The tolerance value was set to $\tau = 10^{-10}$. Hence, the algorithm needs 27 steps to reach this value.

\begin{figure}[h]		
	\hspace{-5mm}	\includegraphics[scale=0.245]{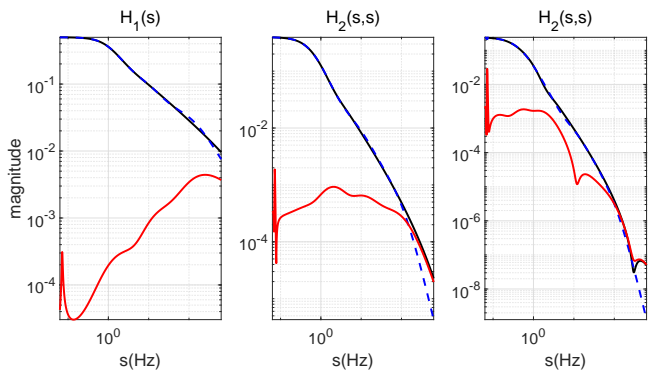}
	\vspace{-5mm}
	\caption{Approximation of the three transfer functions.}
	\label{fig:5}
	\vspace{-3mm}
\end{figure}

In Fig.\;\ref{fig:5}, the errors between the original large-scale system and the reduced one are depicted by evaluating the transfer functions on a denser grid containing 500 points.
%The error increases for the high frequency bandwidths where the level of noise has perturbed a lot the estimations.
\vspace{-1mm}

\begin{figure}[h]
	\hspace{-5mm}		\includegraphics[scale=0.245]{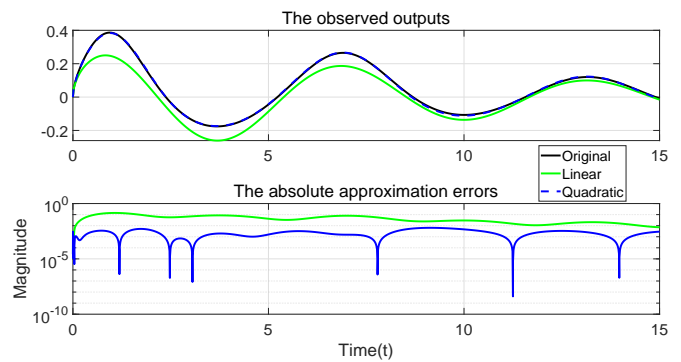}
	\vspace{-5mm}
	\caption{Time-domain simulation: observed outputs (up) and errors (down).}
	\label{fig:6}
	\vspace{-4mm}
\end{figure}

%Finally, in Fig.\;\ref{fig:5} and in Fig.\;\ref{fig:6} qualitative results for both, the frequency domain and the time domain are presented.
Finally, in Fig.\;\ref{fig:6}, we display the time-domain simulations when the control input is given by  $u(t)=\cos(t) e^{-0.1t}$ for $t\in[0,15]$s. Note that the fitted quadratic model improves the approximation accuracy, as compared to the linear one. 

\section{Conclusion}\label{sec:6}

In this paper we have proposed an iterative-based procedure for learning reduced-order quadratic models from data. Instead of employing a one step procedure as was previously proposed in other works, we devised an adaptive scheme that improves the fitted model by including information from higher-order transfer functions. The tests performed show promising results, especially with respect to robustness of the noisy case.  Issues regarding the convergence of the proposed algorithm and of solving ill-conditioned LS problems will be addressed in future works. Possible research endeavors include extending this procedure to other classes of nonlinear control systems, such as quadratic-bilinear or general polynomial systems.

\section{Appendix}\label{sec:7}

\subsection{Proof of Proposition \ref{prop:1} }

Using (\ref{second_trf_quad}), one can write $\bH_2(j \omega)$ as follows
\begin{align} \label{eq1}
\begin{split}
\bH_2(j \omega) &= \bJ_1(2j \omega) \bQ [ \bG_1(j \omega) \otimes \bG_1(j \omega) ].
\end{split}
\end{align}
From (\ref{eq1}), using properties of Kronecker product,we get that
\begin{equation}\label{H2_val_omega}
\left( \bG_1^T(j \omega) \otimes \bG_1^T(j \omega) \otimes \bJ_1(2j \omega) \right) \bv_\bQ = \bH_2(j \omega).
\end{equation}
For all $\omega=\omega_\ell$, with $ 1 \leq \ell \leq N$,  collect the $N$ equalities in (\ref{H2_val_omega}) into matrix format, and hence, Proposition \ref{prop:1} is proven.

\subsection{Proof of Proposition  \ref{prop:2} }

One can write the  function $\bG_2$, evaluated at $j \omega$, as
\begin{align*}
\begin{split}
\bG_2(j &\omega) =  \bPhi(2j\omega) \bQ \Big{[} \bPhi(j\omega) \bB \otimes  \bPhi(j\omega) \bB \Big{]},
\end{split}
\end{align*}
and by employing vectorization techniques, it follows that
\begin{align}\label{eq_G2_1}
\begin{split}
% &\Rightarrow \bG_2(j \omega)  = \left[ \left( \bG_1^T(j \omega) \otimes \bG_1^T(j \omega) \right)  \otimes \bPhi(2j \omega) \right] \bv_{\bQ} \\ 
& \bG_2^T(j \omega)  = \bv_{\bQ}^T \left[ \left( \bG_1(j \omega) \otimes \bG_1(j \omega) \right) \otimes \bPhi^T(2j \omega) \right]. 
\end{split}
\end{align}
One can write the third transfer function  evaluated at $j \omega$
\begin{align*}
\begin{split}
\bH_3(j \omega) &= 2\bJ_1(3j \omega) \bQ \left[\bG_2(j\omega) \otimes \bG_1 (j \omega)\right],
\end{split}
\end{align*}
and by using properties of the Kronecker product, write
\begin{equation} \label{eq_H3_2}
2\left( \bG_2^T(j \omega) \otimes \bG_1^T(j \omega) \otimes \bJ_1(3j \omega) \right) \bv_\bQ = \bH_3(j \omega).
\end{equation}
Substituting (\ref{eq_G2_1}) into (\ref{eq_H3_2}), it follows that
\begin{equation*} 
\bv_{\bQ}^T  \bK(j\omega)  \bv_\bQ = \bH_3(j \omega).
\end{equation*}

\small

\bibliographystyle{spmpsci}  

\bibliography{GKA_ECC_2021_ref}             % bib file to produce the bibliography

\begin{thebibliography}{10}
\providecommand{\url}[1]{{#1}}
\providecommand{\urlprefix}{URL }
\expandafter\ifx\csname urlstyle\endcsname\relax
  \providecommand{\doi}[1]{DOI~\discretionary{}{}{}#1}\else
  \providecommand{\doi}{DOI~\discretionary{}{}{}\begingroup
  \urlstyle{rm}\Url}\fi

\bibitem{ACA05}
Antoulas, A.C.: Approximation of large-scale dynamical systems.
\newblock SIAM, Philadelphia (2005)

\bibitem{ABG20}
Antoulas, A.C., Beattie, C.A., Gugercin, S.: Interpolatory Methods for Model
  Reduction.
\newblock SIAM, Philadelphia (2020)

\bibitem{AGI16}
Antoulas, A.C., Gosea, I.V., Ionita, A.C.: Model reduction of bilinear systems
  in the {L}oewner framework.
\newblock SIAM Journal on Scientific Computing \textbf{38(5)}, B889--B916
  (2016)

\bibitem{ALI17}
Antoulas, A.C., Lefteriu, S., Ionita, A.C.: A tutorial introduction to the
  {L}oewner framework for model reduction.
\newblock In: Model Reduction and Approximation, chap.~8, pp. 335--376. SIAM
  (2017)

\bibitem{astolfi2010}
Astolfi, A.: Model reduction by moment matching for linear and nonlinear
  systems.
\newblock IEEE Transactions on Automatic Control \textbf{55}(10), 2321--2336
  (2010)

\bibitem{BBF14}
Baur, U., Benner, P., Feng, L.: Model order reduction for linear and nonlinear
  systems: A system-theoretic perspective.
\newblock Archives of Computational Methods in Engineering \textbf{21}(4),
  331--358 (2014)

\bibitem{BR71}
Bedrosian, E., Rice, S.: The output properties of {V}olterra systems driven by
  harmonic and {G}aussian inputs.
\newblock Proceedings of the IEEE \textbf{59}, 1688--1707 (1971)

\bibitem{BB15}
Benner, P., Breiten, T.: Two-sided projection methods for nonlinear model order
  reduction.
\newblock SIAM Journal on Scientific Computing \textbf{37}(2), B239--B260
  (2015)

\bibitem{BGG18}
Benner, P., Goyal, P., Gugercin, S.: $\mathcal{H}_2$-quasi-optimal model order
  reduction for quadratic-bilinear control systems.
\newblock SIAM Journal on Matrix Analysis and Applications \textbf{39}(2),
  983--1032 (2018)

\bibitem{morBenGKetal20}
Benner, P., Goyal, P., Kramer, B., Peherstorfer, B., Willcox, K.: Operator
  inference for non-intrusive model reduction of systems with non-polynomial
  nonlinear terms.
\newblock CompMethAppMechEng \textbf{372} (2020)

\bibitem{BOCW17}
Benner, P., Ohlberger, M., Cohen, A., Willcox, K.: Model Reduction and
  Approximation.
\newblock Society for Industrial and Applied Mathematics, Philadelphia, PA
  (2017).
\newblock \doi{10.1137/1.9781611974829}

\bibitem{BoydChua85}
Boyd, S., Chua, L.S.: Fading memory and the problem of approximating nonlinear
  operators with volterra series.
\newblock IEEE Transactions on Circuits and Systems \textbf{30}, 1150--1161
  (1985)

\bibitem{GA18}
Gosea, I.V., Antoulas, A.C.: Data-driven model order reduction of
  quadratic-bilinear systems.
\newblock Numerical Linear Algebra with Applications \textbf{25}(6), e2200
  (2018)

\bibitem{morGosKA20}
Gosea, I.V., Karachalios, D., Antoulas, A.C.: Modeling in the {L}oewner
  framework: from linear dynamics to quadratic nonlinearities.
\newblock In: 21st IFAC World Congress, Berlin, Germany (2020).
\newblock Extended abstract

\bibitem{GS99}
Gustavsen, B., Semlyen, A.: Rational approximation of frequency domain
  responses by vector fitting.
\newblock IEEE Trans. Power Delivery \textbf{14}(3), 1052--1061 (1999)

\bibitem{morKarGA20}
Karachalios, D.S., Gosea, I.V., Antoulas, A.C.: On bilinear time domain
  identification and reduction in the {L}oewner framework.
\newblock In: Model Reduction of Complex Dynamical Systems, International
  Series of Numerical Mathematics. Springer (2020).
\newblock Accepted

\bibitem{MA07}
Mayo, A.J., Antoulas, A.C.: A framework for the solution of the generalized
  realization problem.
\newblock Linear Algebra and Its Applications \textbf{425}(2-3), 634--662
  (2007)

\bibitem{PGW17}
Peherstorfer, B., Gugercin, S., Willcox, K.: Data-driven reduced model
  construction with time-domain {L}oewner models.
\newblock SIAM Journal on Scientific Computing \textbf{39}(5), A2152--A2178
  (2017)

\bibitem{PEHERSTORFER2016196}
Peherstorfer, B., Willcox, K.: Data-driven operator inference for nonintrusive
  projection-based model reduction.
\newblock Computer Methods in Applied Mechanics and Engineering \textbf{306},
  196 -- 215 (2016).
\newblock \doi{https://doi.org/10.1016/j.cma.2016.03.025}

\bibitem{PBK16}
Proctor, J., Brunton, S.L., Kutz, J.N.: Dynamic mode decomposition with
  control.
\newblock SIAM Journal on Applied Dynamical Systems \textbf{15}(1), 142--161
  (2016)

\bibitem{Qin06}
Qin, S.J.: An overview of subspace identification.
\newblock Computers and Chemical Engineering \textbf{30}, 1502--1513 (2006)

\bibitem{WJR05}
Rugh, W.J.: Nonlinear system theory - The Volterra/Wiener approach.
\newblock University Press (1981)

\bibitem{morSch10}
Schmid, P.J.: Dynamic mode decomposition of numerical and experimental data.
\newblock Journal of Fluid Mechanics \textbf{656}, 5--28 (2010)

\bibitem{OM94}
Van~Overschee, P., De~Moor, B.: Subspace algorithms for the identification of
  combined deterministic-stochastic systems.
\newblock Automatica \textbf{30}(1), 75--93 (1994)

\end{thebibliography}

\end{document}